\documentclass{amsart}
\usepackage{fourier}
\usepackage[scaled=0.75]{helvet}

\usepackage{amssymb,amsthm,amsmath}
\usepackage{bm}         
\usepackage{mathrsfs}   
\usepackage{graphicx}
\usepackage{stmaryrd}   
\usepackage{xcolor}
\usepackage[hidelinks]{hyperref}

\def\Z{\mathbb Z}
\def\Q{\mathbb Q}
\def\R{\mathbb R}
\def\A{\{0, 1\}}           
\def\AA{\mathcal{A}}          
\def\TT{\mathcal{T}}          

\newcommand{\abs}[1]{\left|{#1}\right|}
\newcommand{\floor}[1]{\left\lfloor{#1}\right\rfloor}
\newcommand{\ceil}[1]{\left\lceil{#1}\right\rceil}
\newcommand{\round}[1]{\left\llbracket{#1}\right\rrbracket}

\def\Card{\mathrm{Card}}
\DeclareMathOperator{\SL}{SL}
\DeclareMathOperator{\Nut}{Nut}

\newcommand{\BDsim}{\stackrel{\mathrm{BD}}{\sim}}

\theoremstyle{plain}
\newtheorem{thm}{Theorem}

\theoremstyle{definition}
\newtheorem{defn}[thm]{Definition}
\theoremstyle{remark}

\title{Revisiting Ammann bars: new aperiodic tile sets}
\author{Shigeki Akiyama, Tadahisa Hamada, Katsuki Ito}
\date{}

\begin{document}
\begin{abstract}
This note announces new aperiodic tile sets and places them
in the history,
which goes back to
the lines that Ammann drew on his tiles as matching rules. We start
from such a pattern of lines, the Ammann bars, visible in any tiling by the recently discovered aperiodic
monotile known as the Smith Turtle, and we isolate its combinatorial
content in an object we call a Sturmian lattice: three families of
parallel lines in which the gap between consecutive lines takes two
values only, distributed according to a single parameter $\alpha$. We
then reverse the construction and build, for every real quadratic irrational
$\alpha$, a finite set of tiles with matching rules that tiles the
plane and forces every one of its tilings to carry a Sturmian lattice
of slope $\alpha$, up to Galois conjugacy.
These tile sets are aperiodic, and their expansion constants realize a
unit of every real quadratic field. The construction does not rely on
substitution or hierarchical structure, and the resulting tile sets need
not enforce a self-similar substitution structure. A classification of Sturmian lattices and
complete proofs are given in the full version.
\end{abstract}

\maketitle

\section{Introduction}
\label{sec:intro}

A finite set of tiles is \emph{aperiodic} if it tiles the plane, but
only in ways having no translational symmetry. The first aperiodic sets
arose from the study of Wang tiles: Berger constructed one to show that
the tiling problem for Wang tiles is undecidable \cite{Berger}.
R.~M.~Robinson later found another aperiodic set \cite{RobinsonRM}, and
Penrose found one in several versions, among them the kite-and-dart pair
and a pair of rhombs \cite{Penrose}. Socolar and Taylor found a
single hexagonal tile with matching rules, whose rules can be enforced
by shape alone if the tile is allowed to be disconnected
\cite{SocolarTaylor}. Smith, Myers, Kaplan and Goodman-Strauss then found
aperiodic monotiles, single tiles homeomorphic to a disk that tile the
plane but never periodically, a possibility Penrose had left open
\cite{Penrose, SMKGS}.

Robert Ammann was one of the pioneers, and found several of the early
aperiodic sets \cite[\S\S10.2--10.4]{GS}. He was not a
professional mathematician: he had studied only a little calculus and some
programming languages, and he sorted mail in a post office.
For years he declined the invitations of the tiling community. Marjorie
Senechal has recorded a dinner at her home in November 1987, with
N.~G.~de Bruijn as fellow guest. Ammann came, and was asked how he had
discovered his tilings \cite{Senechal}:
\begin{quote}\itshape
I'd been thinking about the lines of red, blue, and yellow dots used to
reproduce color photographs in newspapers. I drew lines crisscrossed at
appropriate angles, and stared at them for awhile. The tiles just popped
out at me.
\end{quote}
The answer is puzzling. Ammann did, however, use such lines as matching rules for
several of his tile sets: one draws segments on each tile, now called
\emph{Ammann bars}, and requires them to join into infinite straight
lines across the edges \cite[\S\S10.3--10.4]{GS}. Read literally,
his answer says that the lines came first, and that the tiles were made
to fit them.

\begin{figure}[!tp]\centering
\includegraphics[width=0.97\linewidth]{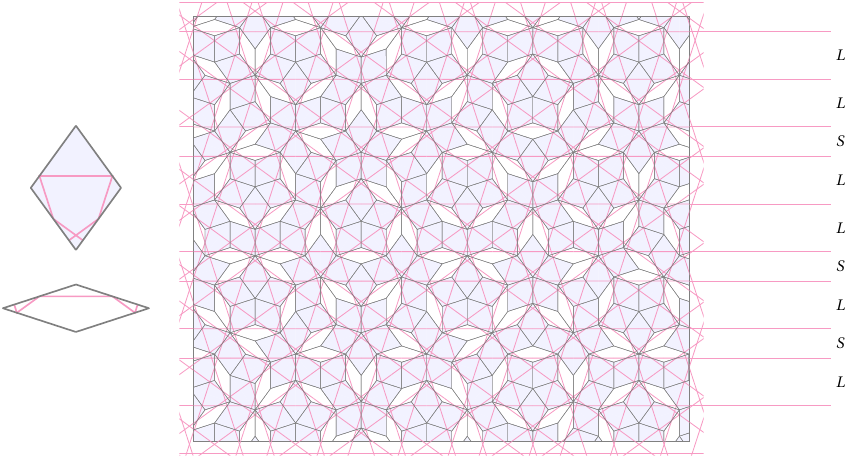}\\[2pt]
{\small (a)}\\[6pt]
\includegraphics[width=0.97\linewidth,page=1]{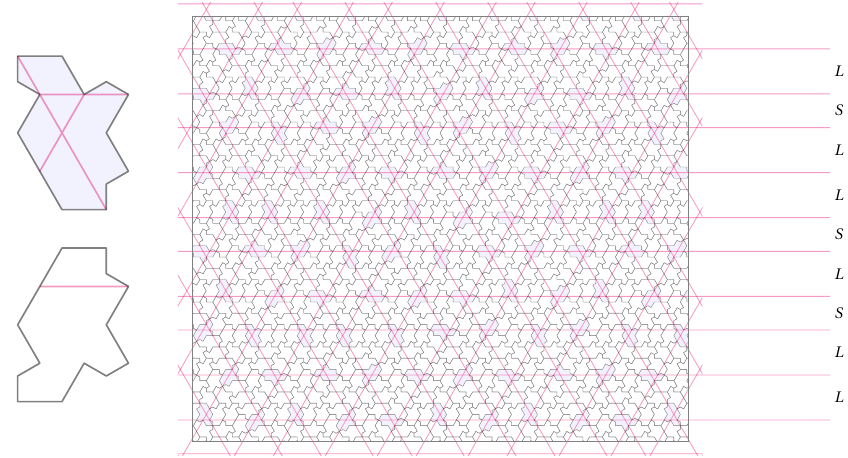}\\[2pt]
{\small (b)}
\caption{Ammann bars (a) on a Penrose tiling by rhombs and
(b) on a Smith Turtle tiling. In (a) the two rhombs on the left carry
the segments of \cite{Fang}; in the tiling every rhomb carries
the same segments, and they join into straight lines running in five
directions. One family is continued to the
right, where its consecutive lines are separated by two distances $L$
and $S$.
In (b) the Turtle (bottom left) carries one segment, and its reflected copy (top left, shaded) carries three, which
bound a small congruent equilateral triangle. In the tiling the segments join into straight lines running in three directions.
The horizontal family is continued to the right; again the gap between consecutive lines takes two values $L$ and
$S$, and so it does in the other two directions.
Whenever two
lines meet, a line of the third direction passes close by, cutting out 
a small equilateral triangle, which lies inside a reflected copy of the Turtle. 
This is the pattern we call a Sturmian
lattice.}
\label{fig:intro}
\end{figure}

Ammann bars can also be drawn on the Penrose tilings \cite{Penrose}.
Gr\"unbaum and Shephard drew them on the kites and darts
\cite[\S10.3]{GS} and left the rhombs as an exercise
\cite[Exercise~10.3.3]{GS}. We draw them on the rhombs with the segments
of Fang, Hammock and Irwin \cite[Fig.~1]{Fang}: every tile of a given
shape carries the same segments, each lying inside the tile, and in any
tiling they
join into straight lines running in five directions
(Figure~\ref{fig:intro}(a)). The lines
of one direction are not equally spaced: the gap between consecutive
lines takes two values only, a long one $L$ and a short one $S$.
Gr\"unbaum and Shephard called the resulting sequences of $L$'s and
$S$'s \emph{musical sequences}, and showed that they are given by the
differences of $\floor{n\tau+c}$, $n\in\Z$, where $\tau=(1+\sqrt{5})/2$
\cite[\S10.6]{GS}; in
present-day language they are Sturmian words. In their view, the
recognition of the importance of Ammann bars was the most significant
advance in the subject. They added \cite[p.~581]{GS}:
\begin{quote}\itshape
it is the systems of bars which are fundamental and the only function
of the tiles is to give a practical realization to them. From this
point of view, the correct approach to the discovery of new aperiodic
sets must be to find tiles which \emph{force} systems of bars with the
required properties.
\end{quote}
Ammann himself used the bars of one of his tile sets to construct an
aperiodic set of sixteen Wang tiles \cite[\S11.1]{GS}.

De~Bruijn's algebraic theory of the Penrose tilings \cite{deBruijn}
was motivated by a model on the line: two-sided sequences of zeros and
ones built by applying a rewriting rule, such as $1\mapsto 10$,
$0\mapsto 1$, over and over again. Like the Penrose tilings, they are
built by repeated substitution, and their finite blocks recur in the
same way as finite patches of a tiling. De~Bruijn noted, however, that
``one thing is missing'' \cite[\S17]{deBruijn}. The matching rules force
the Penrose pattern: any tiling by Penrose's tiles that obeys them is a
Penrose tiling. For the sequences this ``does not seem to have an
analog''; no rule on neighboring symbols characterizes them. Ammann
bars supply the missing analog in the plane. The matching rules force
the bars, and the gaps between the bars of one direction form a musical
sequence, a sequence of de~Bruijn's kind for the rule $L\mapsto LS$,
$S\mapsto L$.

The purpose of this note is to carry out this program literally. We
prescribe a kind of pattern of lines, and we produce a finite set of
tiles, with matching rules, that tiles the plane and forces every one
of its tilings to carry a pattern of that kind. The tile sets obtained
this way are aperiodic, and infinitely many of them are essentially
different.

Which pattern of lines? We take the one hidden in the Smith Turtle,
one of the aperiodic monotiles of \cite{SMKGS}.
Akiyama and Araki gave an alternative proof of its aperiodicity by
drawing Ammann bars on it, in the manner of Ammann \cite{Akiyama-Araki}.
In a Smith Turtle tiling these bars assemble into straight lines running
in three directions (Figure~\ref{fig:intro}(b)). As in the Penrose case,
the gap between consecutive lines of one direction takes two values
only. What is new is the way the three directions fit together: whenever
two lines meet, a line of the third direction passes close by, cutting
out a small equilateral triangle, and these triangles fall on one side
or the other of a given line with no evident rule.

We take the second property as a definition. A \emph{Sturmian lattice}
(\S\ref{sec:SL}) is a family of lines in three directions in which any
three lines that could be expected to be concurrent bound a
triangle of one fixed size instead. This single requirement is already strong.
The gaps in each direction then take at most three values. Apart from
degenerate cases, which are classified explicitly in the full version
and play no role here, there are exactly two gap sizes, which we write as the letters $0$ and $1$. The
three resulting bi-infinite words over $\{0,1\}$ share one common
frequency $\alpha\in[0,1]$ of the letter $1$, the \emph{slope} of the
lattice. For the lattice of slope $\sqrt{6}-2\approx 0.449$ treated in
\S\ref{sec:construction}, for instance, the word of gaps along one
direction contains blocks such as
\[
0\,1\,0\,1\,0\,1\,0\,1\,0\,0\,1\,0\,1\,0\,1\,0\,1\,0\,0\,1,
\]
and in this word the numbers of $1$'s in any two blocks of the same
length differ by at most one. For irrational $\alpha$ these words are Sturmian
words of slope $\alpha$, and the lattice has no nontrivial period.

In the Smith Turtle the tiles come first, and the Sturmian lattice is
read off from their tilings. Our aim is to go the other way: we start
from a Sturmian lattice and \emph{produce} aperiodic tile sets whose
tilings all carry such a lattice. Our main result
(Theorem~\ref{thm:main}) reads as follows.
\begin{quote}\itshape
For every quadratic irrational slope $\alpha$ there is an aperiodic tile
set all of whose tilings carry a Sturmian lattice of slope $\alpha$, up
to Galois conjugacy.
\end{quote}
Why quadratic irrationals? A Sturmian lattice of slope $\alpha$ divides
the plane into cells of three types, with densities proportional to
$(1-\alpha)^{2}:2\alpha(1-\alpha):\alpha^{2}$, a ratio quadratic in
$\alpha$. In our construction, Ammann bars alone do not fix the slope,
and a rational slope would allow periodic tilings. To fix it, we cut a small disk, a
\emph{Bolt}, out of each tile carrying the bars; the Bolts come in
three types, one for each type of cell. We then glue Bolts together
into single tiles; each contains Bolts of the three types in one of
only two prescribed combinations of numbers (\S\ref{sec:construction}). These integer numbers constrain
the three densities, which forces
$\alpha$ to satisfy a quadratic equation with integer coefficients
(\S\ref{sec:construction}). Letting
$\alpha$ range over the quadratic irrationals in $[0,1]$, one obtains
infinitely many essentially different aperiodic tile sets, whose
expansion constants include a unit of every real quadratic field.
Within our current framework, a Sturmian lattice of non-quadratic
irrational slope cannot be enforced \cite[Theorem~4(3)]{Akiyama-Hamada-Ito}. 
However we may have a hope. The restriction comes
from the method, not from Sturmian lattices, which exist for every
irrational slope.

Given a substitutive structure, several methods have been proposed to
construct aperiodic tile sets \cite{Mozes, GoodmanStrauss}; they enforce
a hierarchical unique-composition (substitution) structure within a tile
set. Our route is a different one. Taking Ammann's remark as a guide,
we start from the lines; we then compare the frequencies of the local patterns occurring in a Sturmian
lattice and extract algebraic constraints among them. In the end, the
proof rests on the elementary fact that a quadratic curve and a
line meet in at most two points. Notably, self-similarity is not used in
the proof, and the tile sets need not enforce a self-similar
substitution structure: the
full version gives one whose tiling space has positive topological
entropy \cite[\S10.1]{Akiyama-Hamada-Ito}, whereas primitive invertible
substitution tiling dynamical systems have zero topological entropy
\cite[Theorem~7.16]{RobinsonEA}. The
key tool is the bounded displacement equivalence of Delone sets
(\S\ref{sec:construction}), which we use to match the tiles arising
in a Sturmian lattice.

The present note states the existence theorem and sketches its proof.
Section~\ref{sec:SL} introduces Sturmian lattices, \S\ref{sec:main}
states the main theorem, and \S\ref{sec:construction} describes the
construction and works out one example. The full classification of
Sturmian lattices, the cardinality estimates, and further examples
appear in \cite{Akiyama-Hamada-Ito}.

\section{Sturmian lattices}
\label{sec:SL}

\subsection*{Tiles and tilings}
We recall the standard notions. Following \cite[Definition~5.2]{BG},
which generalizes \cite{GS}, a \emph{tile} is a pair $(T,a)$ of a
compact set $T\subset\R^{2}$ equal to the closure of its interior and a
\emph{color} $a$ in a finite set $F$.
We do not assume any further topological property of tiles; in
particular a tile may be disconnected, and the tiles of
\S\ref{sec:construction} are (Figure~\ref{fig:ex}).
A \emph{prototile set} $\AA$ is a
finite set of tiles; a \emph{patch} is a collection of images of the
tiles under isometries with pairwise disjoint interiors, and a
\emph{tiling} $\TT$ is a patch whose supports cover $\R^{2}$. A vector
$u$ with $\TT+u=\TT$ is a \emph{period} of $\TT$, and $\TT$ is
\emph{non-periodic} when $0$ is its only period. As in \S\ref{sec:intro},
$\AA$ is \emph{aperiodic} if it admits a tiling but only non-periodic
ones. The Ammann bars of \S\ref{sec:intro} are matching rules of this
kind: the tiles are marked with segments that must extend to straight
lines across the edges \cite[\S10.4]{GS};
such markers or decorations are understood as both colors associated with
tiles in the finite set $F$ and corresponding adjacency conditions.

\subsection*{Balanced and Sturmian words}
The gaps between consecutive parallel lines will be encoded as words
in the letters $0$ and $1$, and the following classical notions describe
how evenly the $1$'s are spread in such a word.
Let $\A^{*}$ and $\A^{\Z}$ be the sets of finite and bi-infinite words, respectively.
For a finite word $w$,
we denote by $\abs{w}$ its length and
by $\abs{w}_{1}$ the number of
$1$'s. A \emph{factor} of a word is a block of consecutive letters in it.
Two words $x, y\in \A^{\Z}$ are said to be \emph{mutually $C$-balanced} if $-C\le \abs{u}_{1} - \abs{v}_{1}\le C$ for all respective factors $u$ and $v$ of $x$ and $y$ with $\abs{u} = \abs{v}$.
A word $w\in \A^{\Z}$ is \emph{$C$-balanced} if $w$ is mutually $C$-balanced with itself.
A $C$-balanced word $w\in\A^{\Z}$ has a well-defined
\emph{slope}, the frequency of the letter $1$:
\[
\alpha(w):=\lim_{N\to\infty}\frac{1}{2N}\,\abs{w_{[-N,N)}}_{1}\in[0,1].
\]
A $1$-balanced word $w$ of irrational slope $\alpha$ is called
\emph{Sturmian}. It is known that for some $\rho\in\R$, either
$w_{n}=\floor{(n+1)\alpha+\rho}-\floor{n\alpha+\rho}$ for all $n$ or
$w_{n}=\ceil{(n+1)\alpha+\rho}-\ceil{n\alpha+\rho}$ for all $n$
\cite[Chapter~6]{Fogg}.

\subsection*{Sturmian lattices}
\begin{figure}[htb]\centering
\includegraphics[width=.55\linewidth]{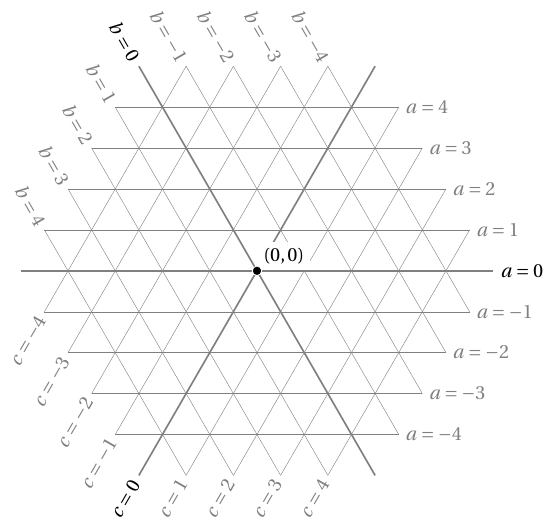}
\caption{Triangular coordinates: the three axes $a,b,c$. The lines
$a=\eta_{0}$, $b=\eta_{1}$, $c=\eta_{2}$ pass through one point when
$\eta_{0}+\eta_{1}+\eta_{2}=0$; otherwise they bound an equilateral
triangle of height $\abs{\eta_{0}+\eta_{1}+\eta_{2}}$, such as the small
ones in Figure~\ref{fig:intro}(b).
}
\label{fig:triplet}
\end{figure}

Consider in $\R^{2}$ three families of parallel lines in the directions
$0$, $2\pi/3$, and $4\pi/3$. Writing
\[
\{a=\eta\}=\{y=\eta\},\quad
\{b=\eta\}=\{y=-\sqrt{3}x-2\eta\},\quad
\{c=\eta\}=\{y=\sqrt{3}x-2\eta\}
\]
for $\eta\in \R$ (Figure~\ref{fig:triplet}),
any three lines $a=\eta_{0}$, $b=\eta_{1}$, $c=\eta_{2}$ bound an
equilateral triangle of height $\abs{\eta_{0}+\eta_{1}+\eta_{2}}$; they
are concurrent precisely when $\eta_{0}+\eta_{1}+\eta_{2}=0$.
Given three bi-infinite real sequences $(a(i))_{i\in\Z}$,
$(b(j))_{j\in\Z}$, $(c(k))_{k\in\Z}$ 
with
\[
a(i+1)-a(i),\quad b(j+1)-b(j),\quad c(k+1)-c(k)\ \ge\ 1,
\]
we form the family of lines
\[
\{a=a(i)\mid i\in\Z\}\cup\{b=b(j)\mid j\in\Z\}\cup\{c=c(k)\mid k\in\Z\}.
\]

\begin{defn}
\label{def:SL}
The family above is a \emph{Sturmian lattice} if 
\[
i+j+k=0\ \Longrightarrow\ \abs{a(i)+b(j)+c(k)}=
\frac{1}{2}.
\]
\end{defn}

Equivalently, every three lines $a=a(i)$, $b=b(j)$, $c=c(k)$ with
$i+j+k=0$ bound a triangle of one fixed size, and these tiny triangles
are arranged, without overlap, in a lattice-like pattern. The Smith Turtle
tiling of Figure~\ref{fig:intro}(b) carries such a structure.

\subsection*{Slope}
Because of the tiny triangles, the parallel lines in a Sturmian lattice
are in general not equally spaced, and the gap between consecutive lines
of one direction takes at most three values.
If it takes exactly one or three values, the lattice is degenerate; such
lattices are described explicitly in \cite[Theorem~3]{Akiyama-Hamada-Ito}
and play no role here.
When there are exactly two, we encode
the two possible gaps between consecutive lines 
in each direction
by letters of $\A$ --- the smaller gap 
by $0$ and the larger
by $1$. This turns $(a(i),b(j),c(k))$ into a triple of
bi-infinite words $(a_{i}),(b_{j}),(c_{k})\in \A^{\Z}$. The defining
condition of Definition~\ref{def:SL} forces these three words to be 
pairwise mutually $1$-balanced, hence to share one common
slope $\alpha\in [0, 1]$, the \emph{slope} of the
Sturmian lattice (the example word in \S\ref{sec:intro} has slope
$\sqrt{6}-2$).
A Sturmian lattice of irrational slope is built from Sturmian
words:

\begin{thm}[{\cite[Theorem~1]{Akiyama-Hamada-Ito}}]
\label{thm:irrational}
Let $(a(i), b(j), c(k))$ be a Sturmian lattice of irrational slope $\alpha$.
Then, there exist $\kappa\ge 1$ and $\bm{\eta} = (\eta_{0}, \eta_{1}, \eta_{2}), \bm{\rho} = (\rho_{0}, \rho_{1}, \rho_{2})\in \R^{3}$ with
\[
\eta_{0} + \eta_{1} + \eta_{2} =
\rho_{0} + \rho_{1} + \rho_{2} = 0
\]
such that
\begin{align*}
a(i)&= i\kappa + \eta_{0} + \round{i\alpha + \rho_{0}},\\
b(j)&= j\kappa + \eta_{1} + \round{j\alpha + \rho_{1}},\\
c(k)&= k\kappa + \eta_{2} + \round{k\alpha + \rho_{2}},
\end{align*}
where $\round{\bullet}\colon \R\to \Z + 1/2$ rounds a real number to the nearest half-integer.
\end{thm}

In words: the lines of each direction sit at the points of an
arithmetic progression of step $\kappa+\alpha$, rounded so that
consecutive gaps are $\kappa$ or $\kappa+1$, the larger gap occurring
with frequency $\alpha$. In particular, any Sturmian lattice of
irrational slope has only the trivial period.

\section{The main theorem}
\label{sec:main}

We now say what it means for tiles to produce a prescribed pattern of
lines. Let $\alpha\in[0,1]$ be a quadratic irrational and $\alpha'$ its
Galois conjugate, the other root of its minimal polynomial. We say that
a prototile set $\AA$ \emph{enforces a Sturmian lattice of slope $\alpha$
up to Galois conjugacy} if $\AA$ admits a tiling and in every tiling by $\AA$ the
Ammann bars form a Sturmian lattice of slope $\alpha$ or $\alpha'$.
(Slopes lie in $[0,1]$, so if $\alpha'\notin[0,1]$ only the slope
$\alpha$ can occur.)

\begin{thm}[{\cite[Theorem~4(2)]{Akiyama-Hamada-Ito}}]
\label{thm:main}
Let $\alpha\in[0,1]$ be a quadratic irrational. Then there exists an
aperiodic tile set $\AA(\alpha)$ that enforces a Sturmian lattice of
slope $\alpha$ up to Galois conjugacy, and whose
underlying expansion constant is a quadratic unit of the real quadratic
field $\Q(\alpha)$.
\end{thm}

We call a tile set $\AA(\alpha)$ as in Theorem~\ref{thm:main} an
\emph{$\SL(\alpha)$-tile set}. A Sturmian lattice of quadratic
irrational slope is asymptotically self-similar, a property inherited
from the self-similarity of Sturmian words
\cite[\S5.3]{Akiyama-Hamada-Ito}, and its \emph{expansion constant} is
the quadratic unit arising from this asymptotic self-similar structure;
for $\alpha=\sqrt{6}-2$ it is $5+2\sqrt{6}$
(\S\ref{sec:construction}). This self-similarity is not used in the
proof, and the tile sets need not enforce a self-similar substitution structure. Since a Sturmian lattice of quadratic
irrational slope has no nontrivial period (\S\ref{sec:SL}), every tiling
by
$\AA(\alpha)$ is non-periodic, so $\AA(\alpha)$ is indeed aperiodic.
As $\alpha$ ranges over the infinitely many quadratic irrationals in
$[0,1]$, Theorem~\ref{thm:main} yields infinitely many aperiodic tile
sets $\AA(\alpha)$, whose expansion constants include a
unit of every real quadratic field.

\section{Construction of the tile sets}
\label{sec:construction}

We outline the construction behind Theorem~\ref{thm:main}; full details
are in \cite{Akiyama-Hamada-Ito}.

The idea is as follows. Ammann bars on the tiles make the lines of a
tiling form a Sturmian lattice, but by themselves they do not fix its
slope. The slope is read off from how often three kinds of local
configurations occur, and it is fixed by gluing these configurations
together, in prescribed numbers, into single tiles. The notion below
says when such a gluing can be done with tiles of bounded size.

\subsection*{Bounded displacement correspondences}
A \emph{Delone set} $X\subset\R^{2}$ is a discrete set whose points are
neither too close nor too sparse: any two are at least some fixed
distance apart, and every disk of some fixed radius contains one. It has a \emph{natural density} $\delta(X)$ whenever the
number of points per unit area converges, which is the case for all the
sets used below. Two Delone sets
$X,Y$ are \emph{bounded displacement equivalent}, written $X\BDsim Y$,
if there is a bijection $\phi\colon X\to Y$ with
$\sup_{x\in X}\abs{\phi(x)-x}<\infty$; then $\delta(X)=\delta(Y)$.

For our purpose, we need a many-to-many version of this concept. Let $X$ and $Y$ be Delone sets with natural densities $\delta(X) : \delta(Y) = p : q$. Suppose that the points of $X$ can be grouped into blocks of $p$ points and those of $Y$ into blocks of $q$ points, and that the blocks can be paired: $X = \bigsqcup_{k\ge 1}X_{k}$ and $Y = \bigsqcup_{k\ge 1}Y_{k}$ with $\Card(X_{k}) = p$ and $\Card(Y_{k}) = q$, such that $\sup_{k}\operatorname{diam}(X_{k}\cup Y_{k})<\infty$. Such partitions play a central role in the construction of our tile sets: each pair $X_{k}\cup Y_{k}$ will become a single tile.
We write $X\xrightarrow{p}\Lambda \xleftarrow{q} Y$, with a dummy Delone set $\Lambda$ indexing the pairs, and call it a \emph{BD correspondence}.
We use the same terminology for the evident generalization to several
Delone sets, grouping prescribed numbers of points from each of them
into uniformly bounded blocks.
In our construction these correspondences are realized by explicit, piecewise-affine maps, inspired by the bounded interpolations between lattices of Duneau and Oguey~\cite{DO}.

\subsection*{Nuts and Bolts}
The prototiles are obtained from a 
hexagon by cutting
out a small disk, which separates it into two parts: the annulus-like
part is a \emph{Nut} and the disk-like part a \emph{Bolt}
(Figure~\ref{fig:NB}). There are three Nuts $\Nut(S)$, $\Nut(M)$, $\Nut(L)$
and three Bolts $S$, $M$, $L$. The three Bolts are disks of different
radii, so that each fits only into the hole of the Nut of its own type.
Nuts and Bolts have complementary roles:
\begin{itemize}
\item the Nuts carry Ammann bars, which govern how Nuts may be placed
side by side; in a valid tiling the bars form a Sturmian lattice;
\item the Bolts are grouped into tiles by a BD correspondence, and this
grouping fixes the slope.
\end{itemize}
In a tiling realizing a Sturmian lattice of slope $\alpha$, the sets
$\mathcal{S}$, $\mathcal{M}$, $\mathcal{L}$ of Bolt centers of each type have
natural densities
\begin{equation}
\label{eq:densities}
\delta(\mathcal{S}):\delta(\mathcal{M}):\delta(\mathcal{L})
=(1-\alpha)^{2}:2\alpha(1-\alpha):\alpha^{2}.
\end{equation}
Thus the proportions of the three kinds of Bolts determine the slope,
and conversely.

\begin{figure}[htb]\centering
\newcommand{\nbpanel}[2]{%
  \begin{minipage}[b]{0.15\linewidth}\centering
    \includegraphics[width=\linewidth,page=#1]{Nuts_and_Bolts-crm.pdf}\\[2pt]
    {\footnotesize #2}
  \end{minipage}}
\nbpanel{1}{$\Nut(S)$}\hfill
\nbpanel{2}{$\Nut(M)$}\hfill
\nbpanel{3}{$\Nut(L)$}\hfill
\nbpanel{4}{Bolt $S$}\hfill
\nbpanel{5}{Bolt $M$}\hfill
\nbpanel{6}{Bolt $L$}
\caption{The three Nuts and the three Bolts. The red segments on the
Nuts are Ammann bars: in a tiling they must continue across the edges
into straight lines, and these lines form a Sturmian lattice. Each Bolt
fills the hole of the Nut of the same type; the dashed outline shows the
Nut around it.}
\label{fig:NB}
\end{figure}

\begin{figure}[htb]\centering
\includegraphics[width=\linewidth,page=1]{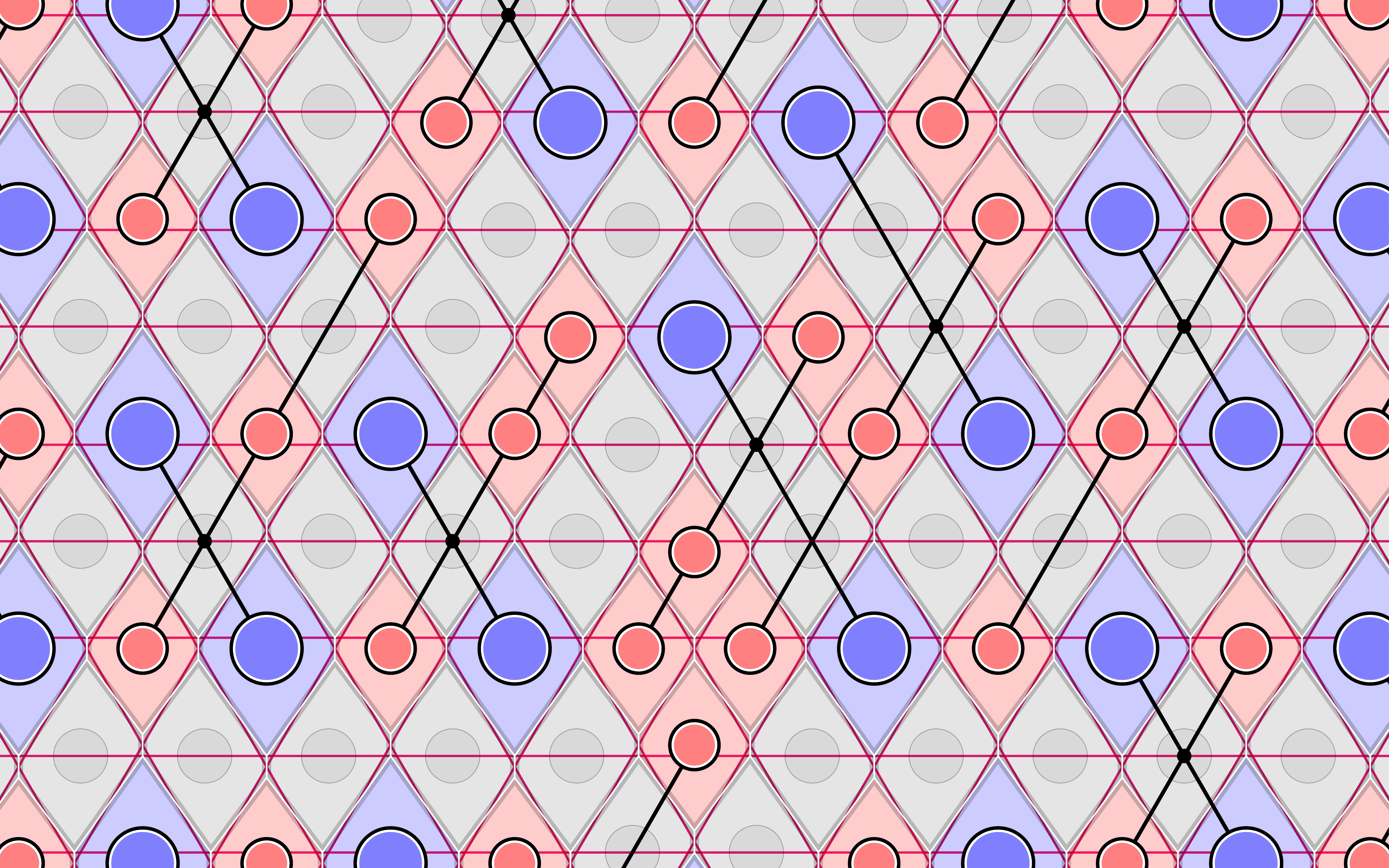}
\caption{Part of a tiling by the $\SL(\sqrt{6}-2)$-tile set of
\S\ref{sec:construction}. The Nuts $\Nut(S)$, $\Nut(M)$, $\Nut(L)$ are
shaded red, gray and blue, and their Ammann bars join into the lines of a
Sturmian lattice of slope $\sqrt{6}-2$. The small red disks are Bolts
$S$, the large blue ones Bolts $L$. The black segments show which Bolts
belong to one tile: each black tree joins three Bolts $S$ and two Bolts
$L$ into a single tile of type $3S+2L$ (Figure~\ref{fig:ex}), while
every Bolt $M$ (gray) is a tile by itself.}
\label{fig:example}
\end{figure}

\subsection*{Construction and aperiodicity}
A \emph{patch-tile} is a single tile made of several Bolts, that is, a
disjoint union of disks (Figure~\ref{fig:ex}).
For the density argument below, only the numbers of Bolts in a
patch-tile matter. We write $xS+yM+zL$ for the class of all patch-tiles made of
$x$ Bolts $S$, $y$ Bolts $M$ and $z$ Bolts $L$, whatever the relative
position of the disks. In terms of Delone sets, replace each Bolt by its
center: the patch-tiles of the class $xS+yM+zL$ in a tiling are then the
blocks of a BD correspondence that groups $x$, $y$, $z$ points of
$\mathcal{S}$, $\mathcal{M}$, $\mathcal{L}$. We therefore call
$xS+yM+zL$ the \emph{BD-equivalence class} of these patch-tiles.
Based on the roles of Nuts and Bolts explained above, we seek a
prototile set $\AA$ consisting of the three Nuts and of patch-tiles of
two classes, that is, of the shape
\begin{equation}
\label{eq:Aform}
\AA/{\BDsim}=\bigl\{\,
T'=x'S+y'M+z'L,\ \ T''=x''S+y''M+z''L,\ \
\Nut(S),\ \Nut(M),\ \Nut(L)\,\bigr\}.
\end{equation}
How should we set $(x', y', z'), (x'', y'', z'')\in \Z_{\ge0}^{3}$ for $\AA$ to enforce a Sturmian lattice of a given quadratic irrational slope $\alpha$ up to Galois conjugacy?

For a valid tiling by $\AA$, every Bolt belongs to a tile of class $T'$
or $T''$, so the density vector
$(\delta(\mathcal{S}),\delta(\mathcal{M}),\delta(\mathcal{L}))$ lies in
the plane spanned by the two integer vectors $(x',y',z')$ and
$(x'',y'',z'')$. By \eqref{eq:densities}, so does the vector
$\bigl((1-\alpha)^{2},2\alpha(1-\alpha),\alpha^{2}\bigr)$,
so it satisfies a linear relation with integer coefficients, not all
zero. Such a relation is a polynomial equation for $\alpha$ of degree at
most two with integer coefficients; this is why, for a tile set of the shape
\eqref{eq:Aform}, an irrational $\alpha$ has to be quadratic.
For every tiling by $\AA$ whose Ammann bars form a Sturmian lattice of
some slope $\beta$, the density vector again lies in this plane, so the
same polynomial vanishes at $\beta$. This polynomial is nonzero, has
degree at most two, and has the quadratic irrational $\alpha$ as a
root, so it is a constant multiple of the minimal polynomial of
$\alpha$; hence $\beta=\alpha$ or $\beta=\alpha'$. Either way $\beta$
is irrational, which is all that aperiodicity requires.
Conversely, for quadratic $\alpha$ the full version chooses suitable
nonnegative integer vectors by an explicit algorithm
\cite[\S8]{Akiyama-Hamada-Ito} and realizes the corresponding partitions
by uniformly bounded BD correspondences, given by explicit
piecewise-affine maps \cite[\S7.3]{Akiyama-Hamada-Ito}; we outline the
latter.

To obtain $\AA$, we prepare Delone sets $\mathcal{S}$, $\mathcal{M}$, and $\mathcal{L}$ on a Sturmian lattice of slope $\alpha$, and partition them.
First, we divide each of them into two subsets $\mathcal{S} = \mathcal{S}'\sqcup \mathcal{S}''$, $\mathcal{M} = \mathcal{M}'\sqcup \mathcal{M}''$, and $\mathcal{L} = \mathcal{L}'\sqcup \mathcal{L}''$ with
\begin{align*}
\delta(\mathcal{S}') :
\delta(\mathcal{M}') :
\delta(\mathcal{L}')&=
x' : y' : z',&
\delta(\mathcal{S}'') :
\delta(\mathcal{M}'') :
\delta(\mathcal{L}'')&=
x'' : y'' : z''.
\end{align*}
Second, we make a BD correspondence between $\mathcal{S}'$, $\mathcal{M}'$, and $\mathcal{L}'$, which groups their points into blocks of $x'$, $y'$, $z'$ Bolts; each block is a patch-tile, and these patch-tiles are the tiles of $\AA$ in the class $T'$. The same is done for $T''$.
Since the blocks are uniformly bounded and Nut-and-Bolt tilings have finite local complexity (for any ball of fixed radius, there are only finitely many patches in the ball up to isometry), only finitely many patch-tiles occur up to isometry, so $\AA$ is finite, as required.

\subsection*{An example: slope $\alpha=\sqrt{6}-2$}
Consider the tile set
\[
\AA/{\BDsim}=\bigl\{\,3S+2L,\ \ M,\ \
\Nut(S),\Nut(M),\Nut(L)\,\bigr\}.
\]
The patch-tile $3S+2L$ forces $\delta(\mathcal{S}):\delta(\mathcal{L})
=3:2$, so we have $(1-\alpha)^{2}:\alpha^{2}=3:2$ by \eqref{eq:densities}, that
is $\alpha^{2}+4\alpha-2=0$ and $\alpha=\sqrt{6}-2$ (the conjugate
$-\sqrt{6}-2$ lies outside $[0,1]$, so here the slope is determined
uniquely); the underlying
expansion constant is the fundamental unit $5+2\sqrt{6}$ of
$\Q(\sqrt{6})$.
Unfolding the BD-equivalence classes gives an
explicit $\SL(\sqrt{6}-2)$-tile set with 29 patch-tiles.
Figure~\ref{fig:ex} shows some tiles of type $3S + 2L$; each is a disjoint union of disks.
Figure~\ref{fig:example} shows how they fit together with the Nuts in a
tiling.

\begin{figure}[htb]\centering
\includegraphics[pagebox=artbox,width=\linewidth,page=7]{Nuts_and_Bolts-crm.pdf}
\caption{A few examples of tiles of type $3S + 2L$. Each tile is the
union of five disjoint disks, three Bolts $S$ and two Bolts $L$; the
dashed outlines indicate the positions of the surrounding Nuts. Tiles of
this type differ in the relative position of the five disks. In the
tiling of
Figure~\ref{fig:example} the five disks of one tile are joined by a
black tree.}
\label{fig:ex}
\end{figure}

The construction subsumes the motivating example. The Smith Turtle
itself corresponds to a quadratic slope in $\Q(\sqrt{5})$
\cite{Akiyama-Araki, Akiyama-Hamada-Ito}, whereas
Theorem~\ref{thm:main} produces an aperiodic tile set for
every quadratic irrational slope.

\section*{Acknowledgements}
The authors are indebted to Professors Michihiro Hirayama and Hajime Kaneko for
their comments on an earlier version of this manuscript, which improved
its readability.
The authors thank the Research Institute for Mathematical Sciences (RIMS),
an International Joint Usage/Research Center located in Kyoto University, for
its support. This research was partially supported by JSPS Grant Numbers
21H00989, 20K03528 and 24K06662, and by JST SPRING Grant Number JPMJSP2124.

\end{document}